\documentclass[a4paper,11pt, leqno]{amsart}

\usepackage[mathscr]{eucal}
\usepackage{amsmath}
\usepackage{amssymb}
\usepackage{amsfonts}
\usepackage{amsthm}

\usepackage{cases}

\usepackage[dvips]{graphics, color}
\theoremstyle{plain}
\newtheorem{theorem}{Theorem}[section]
\newtheorem{corollary}{Corollary}[section]
\newtheorem{lemma}{Lemma}[section]
\newtheorem{remark}{Remark}[section]

\DeclareMathOperator{\sech}{sech}

\newtheorem{definition}{Definition}[section]

\numberwithin{equation}{section}

\usepackage{color}
\usepackage{ulem}

\newcommand{\Add}[1]{\textcolor{red}{#1}}
\newcommand{\Erase}[1]{\textcolor{red}{\sout{\textcolor{black}{#1}}}}

\def\<{\left<} \def\>{\right>}

\def\bea{\begin{eqnarray} }
\def\eea{\end{eqnarray} }

\def\be{\begin{equation} }
\def\ee{\end{equation} }
\def\qed{\ifhmode\unskip\nobreak\fi\ifmmode\ifinner\else\hskip5pt
\fi\fi\hbox{\hskip5 pt \vrule width4 pt height6 pt depth1.5 pt \hskip1pt }}

\begin{document}

\title[]{Tangentially biharmonic Lagrangian $H$-umbilical
submanifolds in complex space forms}
\author[]{Toru Sasahara}
\address{Center for Liberal Arts and Sciences, 
Hachinohe Institute of Technology, 
Hachinohe, Aomori 031-8501, Japan}
\email{sasahara@hi-tech.ac.jp}

\date{}
\footnote[0]{Abh. Math. Semin. Univ. Hamg. (2015) 85:107-123. Some minor mistakes
 are corrected.
A section is added at the end.}

\begin{abstract}
{\footnotesize
The notion of Lagrangian $H$-umbilical submanifolds was  introduced by B. Y.  Chen
in 1997, and these submanifolds 
have appeared in 
several important problems in the study of Lagrangian submanifolds from the Riemannian geometric
 point of view.
Recently, the author introduced the notion of  tangentially biharmonic submanifolds, 
which are defined as submanifolds such that 
the bitension field of 
the inclusion map has vanishing tangential component.
The normal bundle of a round hypersphere in $\mathbb{R}^n$ can be immersed as 
 a tangentially biharmonic Lagrangian $H$-umbilical
submanifold
in $\mathbb{C}^n$.  
Motivated by this fact, we 
 classify tangentially biharmonic Lagrangian $H$-umbilical submanifolds in complex space forms.}
\end{abstract}


 \maketitle



 \section{Introduction}
A biharmonic map is defined as 
a critical point of the bienergy for all variations.    A map is biharmonic if and only if its
bitension field vanishes.
A submanifold is called a biharmonic submanifold
if its inclusion map is a biharmonic map with respect to the induced metric.
Minimal submanifolds are biharmonic. 
Substantial progress has been made toward classifying non-minimal biharmonic submanifolds
in manifolds with special metric properties (e.g., real space forms, complex space forms, Sasakian
space forms, conformally flat spaces, etc.) since 2000. For recent papers on this topic, see, e.g., \cite{bal}, \cite{fo}, \cite{flmo}
and \cite{liang}.

 

On the other hand, the author \cite{sa} introduced the notion of {\it tangentially biharmonic submanifold}
as an  extension of the notion of biharmonic submanifold,
which is a submanifold such that  the bitension field of the inclusion map has vanishing tangential part. 
We would like to mention that 
this notion coincides with the notion
 of  {\it biconservative submanifold}  introduced by Caddeo et al. in \cite{cad}.


It is known that the normal bundle of a submanifold in Euclidean $n$-space $\mathbb{R}^n$ can be immersed as a Lagrangian submanifold
in complex Euclidean $n$-space $\mathbb{C}^n$ (see \cite{hl}).   
Harvey and Lawson \cite{hl} showed that the normal bundle 
$T^{\perp}M^2$ of  a surface $M^2$  in $\mathbb{R}^3$ is minimal in $\mathbb{C}^3$
 if and only if $M^2$ is minimal. As an extension of this result,
 it was proved in \cite{sa} that $T^{\perp}M^2$ is  a tangentially biharmonic Lagrangian submanifold
 in $\mathbb{C}^3$ if and only if $M^2$ is either minimal, 
 a part of a round sphere, or a part of a
circular cylinder. In the two cases when $M^2$ is not minimal, $T^{\perp}M^2$ is (non-biharmonic) tangentially biharmonic.
In particular, it should be noted that,  if $M^2$ is a part of a round sphere,
 then $T^{\perp}M^2$ is a Lagrangian $H$-umbilical submanifold introduced by B. Y. Chen
 in  \cite{ch 3}.  Similar properties hold for a round hypersphere in $\mathbb{R}^n$ with $n>3$.
We remark  that Lagrangian $H$-umbilical submanifolds can be regarded as  the simplest
 Lagrangian submanifolds next to the totally geodesic ones, and these submanifolds appear in
 several important problems in the study of Lagrangian submanifolds from the Riemannian geometric
 point of view.

  Motivated by the above mentioned facts, this paper classifies
   tangentially biharmonic Lagrangian $H$-umbilical
 submanifolds in complex space forms.

\section{Preliminaries}
 Let $M^n$ be  
 an $n$-dimensional  submanifold of a Riemannian manifold $\tilde M$.
 Let us denote by $\nabla$ and
 $\tilde\nabla$ the Levi-Civita connections on $M^n$ and 
 $\tilde M$, respectively. The
 Gauss and Weingarten formulas are respectively given by
\bea 
 \begin{split}
 \tilde \nabla_XY&= \nabla_XY+h(X,Y), \label{gawe}\\
 \tilde\nabla_X \xi&= -A_{\xi}X+D_X\xi 
 \end{split}
\eea 
 for tangent vector fields $X$, $Y$ and normal vector field $\xi$,
 where $h,A$ and $D$ are the second fundamental 
 form, the shape operator and the normal
 connection.
 The second fundamental form $h$ and the shape operator $A$ are related by
 \be
 \<A_{\xi}X, Y\>=\<h(X, Y), \xi\>. \label{relation}
 \ee
 The mean curvature vector field $H$ is defined by 
$H=(1/n){\rm trace} h.$
 The function $|H|$ is called the  {\it mean curvature}.
 If it vanishes identically, then $M^n$ is called a {\it minimal submanifold}.
 In particular, if $h$ vanishes identically,  then $M^n$ is called a {\it totally geodesic submanifold}.

 Let $\tilde M^n(4\epsilon)$ be a complex space form of complex dimension $n$ and constant
 holomorphic sectional curvature $4\epsilon$. 
 The curvature tensor $\tilde R$ of $\tilde M^n(4\epsilon)$ is given by
 \begin{align*}
 \tilde R(X, Y)Z&=\epsilon\{\<Y, Z\>X-\<X, Z\>Y+\<Z, JY\>JX\nonumber\\
 &-\<Z, JX\>JY+2\<X, JY\>JZ\},\label{curvature}\end{align*}
 where $\<, \>$ is the inner product and $J$ is the complex structure of $\tilde M^n(4\epsilon)$.
 
 Any complete and simply connected complex space form $\tilde M^n(4\epsilon)$
 is  holomorphically isometric to the complex Euclidean space $\mathbb{C}^n$, 
 the complex projective space $\mathbb{C}P^n(4\epsilon)$ or the complex hyperbolic space $\mathbb{C}H^n(4\epsilon)$
 according as  $\epsilon=0$, $\epsilon>0$ or $\epsilon<0$. 
 A submanifold $M^n$ of  $\tilde M^n(4\epsilon)$ is called 
 $Lagrangian$ if  $\<X, JY\>=0$ for all tangent vector fields $X$ and $Y$ of $M^n$.
 
 Let $M^n$ be   a Lagrangian submanifold of $\tilde M^n(4\epsilon)$.
 Since $\tilde\nabla J=0$ holds, it follows from (\ref{gawe}) that $M^n$ satisfies
 \be
 D_XJY=J(\nabla_XY).\label{D}
 \ee
 Denote by $R$ and $R^D$  the Riemann curvature tensor of $\nabla$ and $D$ respectively.
 Then the equations of
 Gauss and  Codazzi are respectively given by
 \begin{align}
 \<R(X,Y)Z,W\>&=
 \epsilon(\<X,W\>\<Y,Z\>-\<X,Z\>\<Y,W\>)+\<[A_{JZ},
 A_{JW}](X),Y\>, \label{Gauss}\\
 ({\bar\nabla}_{X}h)(Y,Z)&=
 ({\bar\nabla}_{Y}h)(X,Z),\label{Codazzi}
 \end{align}
 where $X,Y,Z,W$ 
  are vectors tangent  to
 $M^n$,  and
 $\bar\nabla h$ is defined by
 \be ({\bar\nabla}_{X}h)(Y,Z)= D_X h(Y,Z) - h(\nabla_X
 Y,Z) - h(Y,\nabla_X Z). \nonumber
 \ee

 We note that for Lagrangian submanifolds the equation (\ref{Codazzi}) of Codazzi  coincides with 
  the equation of  Ricci.

 \section{Tangentially biharmonic submanifolds}
 Let $f:(M^n, g) \rightarrow N$ be a smooth map between  two Riemannian manifolds.
 The {\it tension field}
 $\tau(f)$ of $f$ is a section
 of the induced vector bundle $f^{*}TN$
 defined by
 $$
\tau(f):=\mathrm{tr}(\nabla^{f} df)=
\sum_{i=1}^{n}\{\nabla^{f}_{e_i}df(e_i)
-df(\nabla_{e_i}e_i)\},
$$
where $\nabla^f$, $\nabla$ and $\{e_i\}$ denote the induced connection, the connection of $M^n$ and
a local orthonormal basis of $M^n$, respectively.
 If $f$ is an isometric immersion,  we have 
 \be\tau(f)=nH.\label{tauH}\ee
 
 The bienergy $E_2(f)$ of $f$ over compact domain $\Omega\subset M^n$ is defined by 
 \be
 E_2(f)=\int_{\Omega}|\tau(f)|^2dv_g,\nonumber
 \ee
 where $dv_g$ is the volume form of $M^n$ (see \cite{es2}).
 If $f$ is a critical point of $E_2$ with respect to compactly supported variations, then $f$ is called a
 {\it biharmonic map} (or {\it $2$-harmonic map}).
Jiang \cite{ji2} proved that $f$ is biharmonic if and only if its bitension field defined by
 \be\tau_2(f):=-\Delta_f\tau(f)+{\rm trace}_gR^{N}(\tau(f),df)df\label{tau2}
 \ee
 vanishes identically, where
 $\Delta_f=-{\rm trace}_g(\nabla^{f}\nabla^{f}-\nabla^{f}_{\nabla})$ and 
 $R^{N}$ is 
 the curvature tensor  of $N$, which is given by $$R^{N}(X, Y)Z=[\nabla^N_X, \nabla^N_Y]Z-\nabla^N_{[X, Y]}Z$$ for the Levi-Civita connection $\nabla^N$ of $N$.
 
 If $f$ is an isometric immersion, then $M^n$ 
 is called a {\it biharmonic submanifold} in $N$.  It follows from (\ref{tauH}) and (\ref{tau2}) that any minimal submanifold is biharmonic. 
 However, the converse is not true in general. 

 In \cite{sa}, the notion of  tangentially biharmonicity, which is weaker than 
  biharmonicity for submanifolds, was introduced as follows:
 \begin{definition}
{\rm  Let $f:M\rightarrow N$ be an isometric immersion. Then $M$
 is called a {\it tangentially biharmonic submanifold} in $N$
 if it satisfies  
 \be
 \{\tau_2(f)\}^{\top}=0,\label{tan1}
 \ee where $\{\cdot\}^\top$ denotes the tangential part of $\{\cdot\}$. }
 \end{definition}

 By Proposition 2.1 in \cite{flmo}, (\ref{tauH}) and (\ref{tan1}), we have
 \begin{lemma} 
 Let $f: M^n\rightarrow \tilde M^{n}(4\epsilon)$ be a Lagrangian isometric immersion. 
 Then $M$ is a tangentially biharmonic submanifold in   $\tilde M^{n}(4\epsilon)$
 if and only if $f$ satisfies
 \be
 4\sum_{i=1}^nA_{D_{e_i}(\tau(f))}e_i+{\rm grad}|\tau(f)|^2=0.\label{TBE}
 \ee
 \end{lemma}

Let $x:M^{n-1}\rightarrow \mathbb{R}^n$ be an isometric immersion.
The normal bundle $T^{\perp}M^{n-1}$ of $M^{n-1}$ is naturally immersed in $\mathbb{R}^n
\times\mathbb{R}^n=\mathbb{R}^{2n}$ by the immersion 
$f(\xi_x):=(x, \xi_x)$, which is  expressed as
\be
f(x, s)=(x, sN)\label{nbundle}
\ee
for the unit normal vector field $N$ along $x$.
We equip $T^{\perp}M^{n-1}$ with the metric induced by $f$. 
If we define the complex structure $J$ on $\mathbb{C}^n=\mathbb{R}^n\times\mathbb{R}^n$
by $J(X, Y):=(-Y, X)$, then 
$T^{\perp}M^{n-1}$ is a Lagrangian submanifold in $\mathbb{C}^n$
(see \cite[III.3.C]{hl}). 
It was proved  in \cite{sa} that  $T^{\perp}M^2$ is a tangentially biharmonic Lagrangian submanifold in
 $\mathbb{C}^3$
 if and only if
$M^2$ is either  minimal, a part of a round sphere or a part of a circular cylinder in $\mathbb{R}^3$.

 
\begin{remark}
{\rm The equation (\ref{TBE}) is equivalent to the condition that a certain stress-energy tensor $S_2$ 
for $E_2(f)$  
is divergence-free.  For a detailed definition of $S_2$, see \cite{cad}.
Caddeo et al. \cite{cad} called these submanifolds satisfying such a condition
{\it biconservative} submanifolds, and moreover they have classified biconservative surfaces in 3-dimensional real space
forms.
 Also, in order to prove the non-existence of non-minimal
  biharmonic hypersurfaces in $\mathbb{E}^4$, Hasanis and Vlachos \cite{hav} classified 
  hypersurfaces 
  satisfying (\ref{TBE}). They called such hypersurfaces {\it $H$-hypersurfaces}.}
\end{remark}

 \begin{remark}
 {\rm The first variation formula of $E_2$ obtained in \cite{ch 3}
 implies that  an isometric immersion $f:M\rightarrow N$ is tangentially biharmonic if and only if it 
 is a critical point of $E_{2}$ 
with respect to all
{\it tangential variations} with compact support.
Here, a tangential variation means a 
variation $f_t$ through $f=f_0$ 
such that the variational vector field 
$V=df_{t}/dt|_{t=0}$ is tangent to $f(M)$.}
 \end{remark}

\begin{remark}
{\rm A submanifold satisfying $\{\tau_2(f)\}^{\perp}=0$ is called a {\it biminimal} submanifold.
Here,  $\{\cdot\}^{\perp}$ denotes the normal part of $\{\cdot\}$ (see \cite{LM3}).}
\end{remark}

 \section{Lagrangian $H$-umbilical submanifolds}
 There exist no totally umbilical Lagrangian submanifolds in a complex space form $\tilde M^n(4\epsilon)$ with
$n\geq 2$ other than the totally geodesic ones (cf. \cite{chogi}).
In view of this fact, Chen \cite{ch 3} introduced the concept of {\it Lagrangian H-umbilical}
submanifolds as the simplest Lagrangian submanifolds next to
 totally geodesic ones in complex space forms as follows.
 
 \begin{definition}
 {\rm A  non-totally geodesic
 Lagrangian  submanifold $M^n$ in a complex space form $\tilde M^n(4\epsilon)$ is called {\it Lagrangian H-umbilical}
 if   
 every point has a neighborhood $W$ on which there 
 exists an  orthonormal frame field $\{e_1, \ldots, e_n \}$
 such that
  the  second fundamental form takes the following form:
 \be
 \begin{split}
 &h(e_1, e_1)=\lambda Je_1, \quad h(e_2, e_2)=\cdots=h(e_n, e_n)=\mu Je_1, \label{lagh}\\
 &h(e_1, e_j)=\mu Je_j, \quad h(e_j, e_k)=0, \quad j\ne k, \quad j, k=2, \ldots, n,
 \end{split}
 \ee 
 where $\lambda$ and $\mu$ are some functions on $W$.

 }
 \end{definition}


Under the assumption that the mean curvature vector field is nowhere vanishing, 
 Lagrangian $H$-umbilical submanifolds  satisfy
 $\hbox{\Erase{$A_H$}}\Add{A_{\frac{H}{|H|}}}JH=\lambda JH$
 for some function $\lambda$ on $M^n$, and moreover, 
 at each point $p$ of  $M^n$,  the shape operator \Erase{$A_{JH}$} \Add{$A_{\frac{H}{|H|}}$}
 has only one eigenvalue $\mu(p)$
 on $D(p)=\{X\in T_pM^n\
  \vert \  \<X, JH\>=0\}$. 
Furthermore, if $\mu$ vanishes nowhere and $\lambda/\mu$ is a constant  $r$ for all points 
$p$ of $M^n$, then $M^n$ is said to {\it be of ratio $r$}. 

 The class of Lagrangian $H$-umbilical submanifolds  includes  important submanifolds. 
 For example, non-minimal twistor holomorphic Lagrangian surfaces in $\mathbb{C}P^2$ (cf. \cite{cas}) and 
  the Whitney's sphere in $\mathbb{C}^n$ (cf. \cite{ch 2}) are Lagrangian $H$-umbilical submanifolds of 
  ratio $r=3$. 
  Due to \cite{chendillen}, a $5$-dimensional Lagrangian submanifold in a complex space
form $\tilde M^5(4\epsilon)$ satisfies the following inequality:
 \be
  \delta(2, 2)\leq (25/4)|H|^2+8\epsilon,\label{delta}\ee
  where $\delta(2, 2)$ is a $\delta$-invariant introduced
   by B.-Y. Chen.  
  A $5$-dimensional Lagrangian $H$-umbilical submanifold  of ratio $4$ 
  in a complex space form $\tilde M^5(4\epsilon)$  satisfies the equality case of (\ref{delta})
  at any point (see \cite{chen1}).

 We need the following definitions to represent Lagrangian $H$-umbilical submanifolds in complex
 space forms.

 \begin{definition}
 {\rm Let  $S^{2n+1}(1)\subset\mathbb{C}^{n+1}$ be  a unit $(2n+1)$-sphere centered at the origin.
 Let $\mathbb{C}^{n+1}_1$ be the 
 complex $(n+1)$-space endowed with the complex coordinates $\{z_1, \ldots, z_{n+1}\}$, whose  inner product
 is given by
$\<(z_1, \ldots , z_{n+1}), (w_1, \ldots, w_{n+1})\>
={\rm Re}(-z_1\bar{w}_1+\sum_{i=2}^{n+1}z_i\bar{w}_i)$.
We put $H_1^{2n+1}(-1)=\{z\in \mathbb{C}^{n+1}_1: \<z, z\>=-1\}$.
 A curve $z=z(s)$ in $S^{3}(1)\subset\mathbb{C}^2$ or
   in $H_1^{3}(-1)\subset\mathbb{C}^2_1$
 is called a {\it Legendre curve} if it satisfies $\<z^{\prime}(s), iz(s)\>=0$ identically.}
\end{definition}
 
 
 \begin{definition}
 {\rm Let $G: N^{n-1}\rightarrow\mathbb{R}^n$ be an isometric immersion of a
 Riemannian $(n-1)$-manifold into Euclidean $n$-space $\mathbb{R}^n$
 and let $F: I\rightarrow\mathbb{C}^{*}$ be a unit speed curve in 
 $\mathbb{C}^{*}:=\mathbb{C}-\{0\}$. Then we can extend the immersion
 $G$ to an immersion  $I\times N^{n-1}$ into ${C}^n$ given by
 \be
F\otimes G: I\times N^{n-1}\rightarrow\mathbb{C}\otimes\mathbb{R}^n=\mathbb{C}^n, \nonumber
 \ee
 where $(F\otimes G)(s, p)=F(s)\otimes G(p)$ for $s\in I$ and $p\in N^{n-1}$.
 We call this extension of $G$ a {\it complex extensor} of $G$ (or of the submanifold $N^{n-1}$)
 via $F$.}
 \end{definition}

 The following theorems are  used for the  classification of  tangentially biharmonic 
 Lagrangian $H$-umbilical submanifolds in complex space forms.

 \begin{theorem}[\cite{ch 3}]\label{t1}
 Let $f: M^n\rightarrow \mathbb{C}P^n(4)$ be a  Lagrangian $H$-umbilical immersion
 whose second fundamental form takes the form  $(\ref{lagh})$, where $n\geq 3$. 
 Assume that the mean curvature vector field is nowhere vanishing. 
 We put $k=e_1\mu/(\lambda-2\mu)$.
 Then $M^n$ satisfies $\mu(\lambda-2\mu)\ne 0$ if and only if up
 to rigid motions of $\mathbb{C}P^n(4)$, $f$ is given by $\pi\circ \psi$, where
 $\pi: S^{2n+1}(1)\subset\mathbb{C}^{n+1}\rightarrow \mathbb{C}P^n(4): 
 w\mapsto w\cdot \mathbb{C}^{*}$ is the Hopf fibration and 
  $\psi$ is given by
 \be
 \psi(s, y_1, \ldots, y_n)=(z_1(s), z_2(s)y_1, \ldots, z_2(s)y_n), \quad y_1^2+\cdots+y_n^2=1,\nonumber
 \ee
 where $z=(z_1, z_2)$ is a unit speed  Legendre curve in $S^3(1)$
 given by
 \begin{align}
&z_1=\frac{i\mu(s)-k(s)}{\sqrt{\mu(s)^2+k(s)^2+1}}
\exp\biggl(\int_0^s(i\lambda(t)-i\mu(t))dt\biggr), \nonumber \\
&z_2= \frac{1}{\sqrt{\mu(s)^2+k(s)^2+1}}\exp\biggl(\int_0^si\mu(t)dt\biggr).
 \nonumber
 \end{align}
 Here $\mu=\mu(s)$ satisfies \be
  \mu^{\prime\prime}(\lambda-2\mu)-\mu^{\prime}(\lambda^{\prime}-3\mu^{\prime})
  +(\lambda-2\mu)^2(-\mu^2+\lambda\mu+1)=0.\label{legen1}
  \ee
 \end{theorem}

 \begin{theorem}[\cite{ch 3}]\label{t2}
 Let $f: M^2\rightarrow \mathbb{C}P^2(4)$ be a Lagrangian $H$-umbilical  immersion whose second fundamental form takes the form  $(\ref{lagh})$. Assume that the mean curvature vector field is nowhere vanishing. 
 Then $M^2$ satisfies $\mu(\lambda-2\mu)\ne 0$ and  the integral curves of $JH$ are geodesics in $M^2$
 if and only if 
 up
  to rigid motions of $\mathbb{C}P^2(4)$, $f$ is given by the immersion described in Theorem {\rm\ref{t1}}  with $n=2$.
 \end{theorem}

 \begin{theorem}[\cite{ch 2}]\label{t3}
 Let $f: M^n\rightarrow \mathbb{C}^n$ be a Lagrangian $H$-umbilical  immersion whose second fundamental form takes the form  $(\ref{lagh})$, where $n\geq 3$. 
 Assume that the mean curvature vector field is nowhere vanishing. 
  Then  $M^n$ satisfies $\mu(\lambda-2\mu)\ne 0$ if and only if  up
  to rigid motions of $\mathbb{C}^n$, $f$ is   a complex extensor of the unit hypersphere in $\mathbb{E}^n$
  via a curve $F=F(s)$ in $\mathbb{C}^{*}$ whose curvature function $\kappa(s)$ and argument $\theta(s)$ satisfy
   $\kappa(s)=\lambda(s)$ and
     $\theta^{\prime}(s)=\mu(s)$. 
 \end{theorem}

 \begin{theorem}[\cite{ch 2}]\label{t4}
  Let $f: M^2\rightarrow \mathbb{C}^2$ be a  Lagrangian $H$-umbilical immersion
  whose second fundamental form takes the form  $(\ref{lagh})$.
  Assume that the mean curvature vector field is nowhere vanishing. 
   Then  $M^n$ satisfies $\mu(\lambda-2\mu)\ne 0$ and the  integral curves of $JH$ are geodesics in $M^2$
   if and only if up
   to rigid motions of $\mathbb{C}^2$, $f$ is a complex extensor of the unit circle in $\mathbb
  {E}^2$
   via  a curve $F=F(s)$ in $\mathbb{C}^{*}$ whose
   curvature function $\kappa(s)$ and argument $\theta(s)$ satisfy
      $\kappa(s)=\lambda(s)$ and
        $\theta^{\prime}(s)=\mu(s)$. 
\end{theorem}

\begin{theorem}[\cite{ch 3}]\label{t5}
 Let $f: M^n\rightarrow \mathbb{C}H^n(-4)$ be a  Lagrangian $H$-umbilical immersion
 whose second fundamental form takes the form  $(\ref{lagh})$, where $n\geq 3$.
 Assume that the mean curvature vector field is nowhere vanishing. 
 We put $k=e_1\mu/(\lambda-2\mu)$.    
Then 
 $M^n$ satisfies $\mu(\lambda-2\mu)\ne 0$ if
 and only if  up
 to rigid motions of $\mathbb{C}H^n(4)$, $f$ is  given by the composition $\pi\circ \psi$, where 
 $\pi: H_1^{2n+1}(-1)\subset \mathbb{C}^{n+1}_1\rightarrow \mathbb{C}H^n(-4):
 w\mapsto w\cdot \mathbb{C}^{*}$ is the Hopf fibration and 
  $\psi$ is one
 of the following immersion:
 
 {\rm (1)} $\mu^2+k^2-1>0$ and 
\be
\psi (s, y_1, \ldots, y_n)=(z_1(s), z_2(s)y_1, \ldots, z_2(s)y_n), \quad y_1^2+y_2^2+\cdots+y_n^2=1,\nonumber
 \ee
 where $z=(z_1, z_2)$ is a unit speed Legendre curve in $H_1^3(-1)$ given by
  \begin{align}
&z_1=\frac{i\mu(s)-k(s)}{\sqrt{\mu(s)^2+k(s)^2-1}}
\exp\biggl(\int_0^s(i\lambda(t)-i\mu(t))dt\biggr),\nonumber\\
&z_2= \frac{1}{\sqrt{\mu(s)^2+k(s)^2-1}}\exp\biggl(\int_0^si\mu(t)dt\biggr).
 \nonumber
 \end{align}
  
 {\rm (2)} $\mu^2+k^2-1<0$ and 
  \be
  \psi(s, y_1, \ldots, y_n)=(z_1(s)y_1, \ldots,  z_1(s)y_n, z_2(s)),  \quad y_1^2-y_2^2-\cdots-y_n^2=1,\nonumber
  \ee
  where $z=(z_1, z_2)$ is a unit speed Legendre curve in $H_1^3(-1)$ given by
   \begin{align}
&z_1=\frac{i\mu(s)-k(s)}{\sqrt{1-\mu(s)^2-k(s)^2}}
\exp\biggl(\int_0^s(i\lambda(t)-i\mu(t))dt\biggr),\nonumber\\
&z_2= \frac{1}{\sqrt{1-\mu(s)^2-k(s)^2}}\exp\biggl(\int_0^si\mu(t)dt\biggr).\nonumber
 \end{align}

 {\rm(3)} $\mu^2+k^2-1=0$ and 
 \be
  \begin{split}
&\psi(s, u_2, \ldots, u_n)\\
=&\exp\biggl[\int_0^s(k(t)+i\mu(t))dt\biggr]\Biggl( 1+\dfrac{1}{2}\sum_{j=2}^nu_j^2-
  \int_0^s(k(x)+i\mu(x)){\rm exp}\biggl(-\int_0^x2k(t)dt\biggr)dx,\nonumber\\
 & [-k(0)+i\mu(0)]\biggl[\dfrac{1}{2}\sum_{j=2}^nu_j^2-\int_0^s(k(x)+i\mu(x))
 \exp\biggl(-\int_0^x2k(t)dt\biggr)dx\biggr], 
  u_2, \ldots, u_n\Biggr).
 \end{split} 
  \ee

  In the case of {\rm (1)} and {\rm (2)}, $\mu=\mu(s)$ satisfies
  \be
  \mu^{\prime\prime}(\lambda-2\mu)-\mu^{\prime}(\lambda^{\prime}-3\mu^{\prime})
  +(\lambda-2\mu)^2(-\mu^2+\lambda\mu-1)=0.\label{legen2}
  \ee
 \end{theorem}

 \begin{theorem}[\cite{ch 3}]\label{t6}
 Let $f: M^2\rightarrow \mathbb{C}H^2(-4)$ be a Lagrangian $H$-umbilical  immersion
 whose second fundamental form takes the form  $(\ref{lagh})$.
 Assume that the mean curvature vector field is nowhere vanishing. 
 Then $M^2$ satisfy $\mu(\lambda-2\mu)\ne 0$ and the integral curves of $JH$ are geodesics in $M^2$ if and 
 only if
up
 to rigid motions of $\mathbb{C}H^2(-4)$,  $f$ is given by 
 one of the immersions described in Theorem {\rm\ref{t5}} with $n=2$.
 \end{theorem}

 \begin{remark}\label{remark1}
 {\rm In the above theorems, $e_1=JH/|H|={\partial}/{\partial s}$.}
 \end{remark}

 \begin{remark}\label{remark2}
 {\rm For a unit speed curve $F=F(s)$ in $\mathbb{C}^{*}$, we put $\alpha=|F|^2$. Then,
\Add{after replacing $s$ by $-s$ if necessary,}
the curvature $\kappa(s)$ and the argument $\theta(s)$ of $F(s)$ is given by (see \cite[p.176]{ch 4})
 \be
 \kappa=\frac{2-\alpha^{\prime\prime}}{\sqrt{4\alpha-\alpha^{\prime 2}}}, \quad
 \theta=\int\frac{\sqrt{4\alpha-\alpha^{\prime 2}}}{2\alpha}ds. \label{kappa}
 \ee
 By a straightforward computation, from (\ref{kappa}) we obtain
 \be\theta^{\prime\prime}=(\kappa-2\theta^{\prime})(\ln |F|)^{\prime}.\label{theta2}\ee
  This implies that
  if $\theta^{\prime}$ is constant
  and $\kappa\ne 2\theta^{\prime}$ , then  $F$ is a circle centered at the origin.
  In this case, it follows from (\ref{kappa}) that $\kappa=\theta^{\prime}=|F|^{-1}$.}
 \end{remark}

 \begin{remark}
 {\rm Equation (\ref{legen1}) (resp. (\ref{legen2}))
  is the necessary and sufficient condition for $z$ in Theorem \ref{t1} (resp. in (1) and (2) of  Theorem \ref{t5})
   to be a Legendre curve of curvature $\lambda(s)$
    in $S^3(1)$ (resp. in $H_1^3(-1)$).
    These equations coincide with the equation of Gauss of $M^n$. Even if two Legendre curves
    are congruent to each other in $S^3(1)$ (resp. in $H_1^3(-1)$),  two Lagrangian $H$-umbilical submanifolds constructed from
      them as in Theorem \ref{t1} (resp. in (1) and (2) of Theorem \ref{t5})
    are not necessarily  congruent to each other.}
 \end{remark}

 \section{Tangentially biharmonic Lagrangian $H$-umbilical submanifolds}
 

Let us denote by $\iota_r$
 the inclusion of a round hypersphere  $S^{n-1}(r)$ of radius $r$ centered at the origin.
 If $x=\iota_r$ in $(\ref{nbundle})$, then 
$T^{\perp}S^{n-1}(r)$ is expressed as 
 a complex extensor of $S^{n-1}(1)$ via $F(s)=r+si$, which is a Lagrangian $H$-umbilical
submanifold of ratio $0$ (see \cite[Example 2.4]{ch 2}). Moreover, by a straightforward computation
we find that it  is a
(non-biharmonic) tangentially biharmonic submanifold in $\mathbb{C}^n$  (see \cite{sa} for $n=3$).
 This motivates us to classify tangentially biharmonic Lagrangian $H$-umbilical submanifolds
in complex space forms.
 
 The main results of this paper are the following  classification theorems.

 \begin{theorem}\label{ct1}
  
  Let  $M^n$ be  a Lagrangian $H$-umbilical submanifold of $\mathbb{C}P^n(4)$ with nowhere
  vanishing mean curvature vector field.
  Then $M^n$ is  tangentially biharmonic if and only if up to rigid motions of 
  $\mathbb{C}P^n(4)$, $M^n$ is locally given by
one of the following:
  
 {\rm (1)} A Lagrangian $H$-umbilical submanifold of  
 constant mean curvature defined by $\pi\circ\psi$,
 where
 \be
  \psi(s, y_1, \ldots, y_n)=
  (z_1(s), z_2(s)y_1, \ldots, z_2(s)y_n), \quad y_1^2+\cdots+y_n^2=1\nonumber
  \ee
 and $z=(z_1, z_2)$ is a unit speed Legendre curve in $S^{3}(1)\subset\mathbb{C}^2$ given by
 \be
 \displaystyle z(s)=\biggl(\frac{i\mu}{\sqrt{\mu^2+1}}e^{-\frac{i}{\mu}s}, \frac{1}{\sqrt{\mu^2+1}}e^{i\mu s}\biggr), \quad \mu\in\mathbb{R}-\{0\}.\nonumber\ee
  In this case, the submanifold has ratio $(1-\mu^{-2})$.

  {\rm (2)} A Lagrangian $H$-umbilical submanifold of ratio $0$ and non-constant
   mean curvature defined by $\pi\circ\psi$,
 where
 \be
  \psi(s, y_1, \ldots, y_n)=
  (z_1(s), z_2(s)y_1, \ldots, z_2(s)y_n), \quad y_1^2+\cdots+y_n^2=1\nonumber
  \ee
 and $z=(z_1, z_2)$ is a unit speed Legendre curve in $S^{3}(1)\subset\mathbb{C}^2$ given by
  \begin{align}
  z_1(s)&=-\frac{2i\mu^2+\mu^{\prime}}{\sqrt{c}\mu^{\frac{3}{2}}}
     \exp\biggl(\int_0^s-i\mu dt\biggr),\nonumber\\
   z_2(s)&=-\frac{2}{\sqrt{c\mu}}\exp\biggl(\int_0^si\mu dt\biggr),\nonumber   
   \end{align}
   with a non-constant positive solution $\mu=\mu(s)$ of
   \be
   \mu^{\prime 2}=-4\mu^2(\mu^2+1)+c\mu^3\nonumber
   \ee
  for some positive constant $c$.

 {\rm  (3)}  A Lagrangian $H$-umbilical submanifold of ratio $(7-n)/3$ and  non-constant
   mean curvature  defined by $\pi\circ\psi$,
 where $n\ne 7$ and 
 \be
  \psi(s, y_1, \ldots, y_n)=
  (z_1(s), z_2(s)y_1, \ldots, z_2(s)y_n), \quad y_1^2+\cdots+y_n^2=1\nonumber
  \ee
 and  $z=(z_1, z_2)$ is a unit speed Legendre curve in $S^{3}(1)\subset\mathbb{C}^2$ given by
 \begin{align}
 z_1(s)&=\frac{i(1-n)\mu^2-3\mu^{\prime}}{3\sqrt{c}\mu^{\frac{n+2}{n-1}}}
 \exp\biggl(\int_0^si\frac{4-n}{3}\mu dt\biggr),\nonumber\\
 z_2(s)&=\frac{(1-n)}{3\sqrt{c}\mu^{\frac{3}{n-1}}}\exp\biggl(\int_0^si\mu dt\biggr),\nonumber
 \end{align}
 with a non-constant positive solution $\mu=\mu(s)$ of
   \be
   \mu^{\prime 2}=-\frac{(1-n)^2}{9}\mu^2(\mu^2+1)+c\mu^{\frac{2(n+2)}{n-1}}\nonumber\ee
   for some  positive constant $c$.

  {\rm (4)} A  Lagrangian $H$-umbilical submanifold of non-constant
   mean curvature defined by $\pi\circ\psi$,
 where
 \be
  \psi(s, y_1, \ldots, y_n)=
  (z_1(s), z_2(s)y_1, \ldots, z_2(s)y_n), \quad y_1^2+\cdots+y_n^2=1\nonumber
  \ee
 and $z=(z_1, z_2)$ is a unit speed Legendre curve in $S^{3}(1)\subset\mathbb{C}^2$ given by
   \begin{align}
&z_1=\frac{i\mu-k}{\sqrt{\mu^2+k^2+1}}
\exp\biggl(\int_0^s(i\lambda-i\mu)dt\biggr),\nonumber\\
&z_2= \frac{1}{\sqrt{\mu^2+k^2+1}}\exp\biggl(\int_0^si\mu dt\biggr).
 \nonumber
 \end{align}
Here $k(s)=\mu^{\prime}(s)/(\lambda(s)-2\mu(s))$, $\lambda\ne r\mu$ for any $r\in \mathbb{R}$, and moreover
  $\lambda=\lambda(s)$ and $\mu=\mu(s)$ are non-constant solutions of 
  \begin{eqnarray}
  \begin{cases}
  \mu^{\prime\prime}(\lambda-2\mu)-\mu^{\prime}(\lambda^{\prime}-3\mu^{\prime})
   +(\lambda-2\mu)^2(-\mu^2+\lambda\mu+1)=0, \\
   (3\lambda+(n-1)\mu)(\lambda-2\mu)\lambda^{\prime}+(n-1)\lambda
   (3\lambda+(n-5)\mu)\mu^{\prime}=0.
  \end{cases}\nonumber
  \end{eqnarray}
  \end{theorem}

 Let $\iota: S^{n-1}(1)\rightarrow\mathbb{E}^n$ be the inclusion of the unit hypersphere
 centered at the origin. 
 For a unit speed curve $F=F(s)$ in $\mathbb{C}^{*}$, let
  $\kappa(s)$ and $\theta(s)$ be the  curvature function
 and the argument of $F$, respectively. Then we have 
 
 \begin{theorem}\label{ct2}
 Let  $M^n$ be  a Lagrangian $H$-umbilical submanifold of $\mathbb{C}^n$ with nowhere
  vanishing mean curvature vector field.
  Then $M^n$ is  tangentially biharmonic if and only if up to rigid motions of 
  $\mathbb{C}^n$, $M^n$ is locally given by
one of the following:

{\rm (1)} A Lagrangian cylinder over a circle:
\be
(a\exp(is/a), u_2, \ldots, u_n), \quad a>0.\nonumber
\ee

 {\rm (2)} A complex extensor $F\otimes\iota$, where $F$ is  a circle centered at the origin.
  
{\rm (3)} A complex extensor $F\otimes\iota$, $F$ is  a  line which does not pass
  through the origin.
  
{\rm (4)} A complex extensor $F\otimes\iota$, 
where  $F$ is a unit speed 
 curve satisfying \be\kappa(s)=((7-n)/3)\theta^{\prime}(s),\quad n\ne \hbox{\Erase{$4$}}, 7, \quad
 \theta^{\prime\prime}(s)\ne 0.  \nonumber\ee
   
{\rm (5)} A complex extensor $F\otimes\iota$, where  $F$ is a unit speed 
 curve satisfying $\kappa\ne r\theta^{\prime}$ for any $r\in \mathbb{R}$ and  \be\kappa^{\prime}(3\kappa+(n-1)\theta^{\prime})+(n-1)({\ln}\hskip2pt |F|)^{\prime}\kappa
   (3\kappa+(n-5)\theta^{\prime})=0.\label{kappatheta}\ee
 
 \noindent Submanifolds of types $(2)$, $(3)$ and 
 $(4)$  are of ratio $1$, $0$ and  $(7-n)/3$, respectively.
 \end{theorem}

 \begin{theorem}\label{ct3}
  Let  $M^n$ be  a Lagrangian $H$-umbilical submanifold of $\mathbb{C}H^n(-4)$
  with nowhere
  vanishing mean curvature vector field.
  Then $M^n$ is tangentially biharmonic if and only if up to rigid motions of 
  $\mathbb{C}H^n(-4)$, $M^n$ is locally given by
one of the following:

{\rm (1)} A flat Lagrangian $H$-umbilical submanifold of ratio $2$ and  constant mean 
 curvature defined by $\pi\circ\psi$, where \be
  \psi(s, u_2, \ldots, u_n)=\dfrac{e^{is}}{\hbox{\Erase{$2$}}}\biggl(1-is+\dfrac{1}{2}\sum_{j=2}^nu_j^2, s+\dfrac{i}{2}\sum_{j=2}^nu_j^2,
  u_2, \ldots, \hbox{\Erase{$n_n$}}\Add{u_n}\biggr).\nonumber
  \ee
  
  {\rm (2)}  A Lagrangian $H$-umbilical submanifold of constant mean curvature
   defined by $\pi\circ\phi$, where 
   \be\psi (s, y_1, \ldots, y_n)=(z_1(s), z_2(s)y_1, \ldots, z_2(s)y_n), \quad y_1^2+y_2^2+\cdots+y_n^2=1\nonumber\ee
 and $z=(z_1, z_2)$ is a unit speed Legendre curve in $H_1^3(-1)$  given by 
 \be
 \displaystyle z(s)=\biggl(\frac{i\mu}{\sqrt{\mu^2-1}}e^{\frac{i}{\mu}s}, \frac{1}{\sqrt{\mu^2-1}}e^{i\mu s}\biggr), \quad \mu^2-1>0, \quad \mu\in\mathbb{R}-\{0\}.\nonumber\ee  
 In this case, the submanifold has ratio $(1+u^{-2})$.

   {\rm (3)} A  Lagrangian $H$-umbilical submanifold of ratio $0$ and 
   non-constant mean curvature 
   defined by $\pi\circ\phi$, where 
   \be\psi (s, y_1, \ldots, y_n)=(z_1(s), z_2(s)y_1, \ldots, z_2(s)y_n), \quad y_1^2+y_2^2+\cdots+y_n^2=1\nonumber\ee
 and $z=(z_1, z_2)$ is a unit speed Legendre curve in $H_1^3(-1)$  given by 
 \begin{align}
  z_1(s)&=-\frac{2i\mu^2+\mu^{\prime}}{\sqrt{c}\mu^{\frac{3}{2}}}
     \exp\biggl(\int_0^s-i\mu dt\biggr),\nonumber\\
   z_2(s)&=-\frac{2}{\sqrt{c\mu}}\exp\biggl(\int_0^si\mu dt\biggr),\nonumber   
   \end{align}
   with a non-constant positive solution $\mu=\mu(s)$ of
   \be
   \mu^{\prime 2}=-4\mu^2(\mu^2-1)+c\mu^3\nonumber
   \ee
  for some  positive constant $c$.


  {\rm  (4)} A Lagrangian $H$-umbilical submanifold of ratio $(7-n)/3$ and 
   non-constant mean curvature
   defined by $\pi\circ\phi$, where $n\ne 7$ and 
   \be\psi (s, y_1, \ldots, y_n)=(z_1(s), z_2(s)y_1, \ldots, z_2(s)y_n), \quad y_1^2+y_2^2+\cdots+y_n^2=1\nonumber\ee
 and $z=(z_1, z_2)$ is a unit speed Legendre curve in $H_1^3(-1)$  given by 
 \begin{align}
 z_1(s)&=\frac{i(1-n)\mu^2-3\mu^{\prime}}{3\sqrt{c}\mu^{\frac{n+2}{n-1}}}
 \exp\biggl(\int_0^si\frac{4-n}{3}\mu dt\biggr),\nonumber\\
 z_2(s)&=\frac{(1-n)}{3\sqrt{c}\mu^{\frac{3}{n-1}}}\exp\biggl(\int_0^si\mu dt\biggr),\nonumber
 \end{align}
 with a non-constant positive solution $\mu=\mu(s)$ of
   \be
   \mu^{\prime 2}=-\frac{(1-n)^2}{9}\mu^2(\mu^2-1)+c\mu^{\frac{2(n+2)}{n-1}}\nonumber\ee
   for some  positive constant $c$.


  
  {\rm (5)}  A Lagrangian $H$-umbilical  submanifold of  
   non-constant mean curvature 
   defined by $\pi\circ\psi$, where 
   \be\psi (s, y_1, \ldots, y_n)=(z_1(s), z_2(s)y_1, \ldots, z_2(s)y_n), \quad y_1^2+y_2^2+\cdots+y_n^2=1\nonumber\ee
 and  $z=(z_1, z_2)$ is a unit speed Legendre curve in $H_1^3(-1)$  given by   
   \begin{align}
&z_1=\frac{i\mu-k}{\sqrt{\mu^2+k^2-1}}
\exp\biggl(\int_0^s(i\lambda-i\mu)dt\biggr),\nonumber\\
&z_2= \frac{1}{\sqrt{\mu^2+k^2-1}}\exp\biggl(\int_0^si\mu dt\biggr).
 \nonumber
 \end{align}
  Here $k(s)=\mu^{\prime}(s)/(\lambda(s)-2\mu(s))$, $\mu^2+k^2-1>0$, 
  $\lambda\ne r\mu$ for any $r\in \mathbb{R}$, and moreover
    $\lambda=\lambda(s)$ and $\mu=\mu(s)$ are non-constant solutions of 
  \begin{eqnarray}
    \begin{cases}
    \mu^{\prime\prime}(\lambda-2\mu)-\mu^{\prime}(\lambda^{\prime}-3\mu^{\prime})
     +(\lambda-2\mu)^2(-\mu^2+\lambda\mu-1)=0, \\
     (3\lambda+(n-1)\mu)(\lambda-2\mu)\lambda^{\prime}+(n-1)\lambda
     (3\lambda+(n-5)\mu)\mu^{\prime}=0.
    \end{cases}\nonumber
    \end{eqnarray}

 {\rm (6)} A Lagrangian $H$-umbilical submanifold of constant mean 
 curvature 
 defined by $\pi\circ\psi$, where
  \be\psi(s, y_1, \ldots, y_n)=(z_1(s)y_1, \ldots,  z_1(s)y_n, z_2(s)),  \quad y_1^2-y_2^2-\cdots-y_n^2=1\nonumber\ee
   and $z=(z_1, z_2)$ is a unit speed Legendre curve in $H_1^3(-1)$ given by 
    \be
    \displaystyle z(s)=\biggl(\frac{1}{\sqrt{1-\mu^2}}e^{i\mu s}, \frac{i\mu}{\sqrt{1-\mu^2}}e^{\frac{i}{\mu}s}\biggr), \quad 1-\mu^2>0, \quad \mu\in\mathbb{R}-\{0\}.\nonumber\ee
    In this case, the submanifold has ratio $(1+\mu^{-2})$. 
     
     {\rm (7)} A Lagrangian $H$-umbilical submanifold of ratio $0$ and non-constant mean 
 curvature 
 defined by $\pi\circ\psi$, where
  \be\psi(s, y_1, \ldots, y_n)=(z_1(s)y_1, \ldots,  z_1(s)y_n, z_2(s)),  \quad y_1^2-y_2^2-\cdots-y_n^2=1\nonumber\ee
  and $z=(z_1, z_2)$ is a unit speed Legendre curve in $H_1^3(-1)$ given by 
  \begin{align}
  z_1(s)&=-\frac{2}{\sqrt{-c\mu}}\exp\biggl(\int_0^si\mu dt\biggr),\nonumber \\  
  z_2(s)&=-\frac{2i\mu^2+\mu^{\prime}}{\sqrt{-c}\mu^{\frac{3}{2}}}
     \exp\biggl(\int_0^s-i\mu dt\biggr),\nonumber
   \end{align}
   with a non-constant positive solution $\mu=\mu(s)$ of
   \be
   \mu^{\prime 2}=-4\mu^2(\mu^2-1)+c\mu^3\nonumber
   \ee
  for some negative constant $c$.


    {\rm  (8)}  A Lagrangian $H$-umbilical submanifolds of ratio $(7-n)/3$ and non-constant mean 
 curvature
 defined by $\pi\circ\psi$, where $n\ne 7$ and 
  \be\psi(s, y_1, \ldots, y_n)=(z_1(s)y_1, \ldots,  z_1(s)y_n, z_2(s)),  
  \quad y_1^2-y_2^2-\cdots-y_n^2=1\nonumber\ee
   and $z=(z_1, z_2)$ is a unit speed Legendre curve in $H_1^3(-1)$ given by
  \begin{align}
  z_2(s)&=\frac{(1-n)}{3\sqrt{-c}\mu^{\frac{3}{n-1}}}\exp\biggl(\int_0^si\mu dt\biggr),\nonumber\\
 z_1(s)&=\frac{i(1-n)\mu^2-3\mu^{\prime}}{3\sqrt{-c}\mu^{\frac{n+2}{n-1}}}
 \exp\biggl(\int_0^si\frac{4-n}{3}\mu dt\biggr),\nonumber 
 \end{align}
 with a non-constant positive solution $\mu=\mu(s)$ of
   \be
   \mu^{\prime 2}=-\frac{(1-n)^2}{9}\mu^2(\mu^2-1)+c\mu^{\frac{2(n+2)}{n-1}}\nonumber\ee
   for some  negative constant $c$.


    
    {\rm (9)} A Lagrangian $H$-umbilical submanifold of  non-constant mean 
 curvature 
 defined by $\pi\circ\psi$, where
  \be\psi(s, y_1, \ldots, y_n)=(z_1(s)y_1, \ldots,  z_1(s)y_n, z_2(s)),  
  \quad y_1^2-y_2^2-\cdots-y_n^2=1\nonumber\ee
   and $z=(z_1, z_2)$ is a unit speed Legendre curve in $H_1^3(-1)$ given by
    \begin{align}
&z_1=\frac{i\mu-k}{\sqrt{1-\mu^2-k^2}}
\exp\biggl(\int_0^s(i\lambda-i\mu)dt\biggr),\nonumber\\
&z_2= \frac{1}{\sqrt{1-\mu^2-k^2}}\exp\biggl(\int_0^si\mu dt\biggr).
 \nonumber
 \end{align} 
    Here $k(s)=\mu^{\prime}(s)/(\lambda(s)-2\mu(s))$, $1-\mu^2-k^2>0$, 
    $\lambda\ne r\mu$ for any $r\in \mathbb{R}$, and moreover
        $\lambda=\lambda(s)$ and $\mu=\mu(s)$ are non-constant solutions of 
      \begin{eqnarray}
        \begin{cases}
        \mu^{\prime\prime}(\lambda-2\mu)-\mu^{\prime}(\lambda^{\prime}-3\mu^{\prime})
         +(\lambda-2\mu)^2(-\mu^2+\lambda\mu-1)=0, \\
         (3\lambda+(n-1)\mu)(\lambda-2\mu)\lambda^{\prime}+(n-1)\lambda
         (3\lambda+(n-5)\mu)\mu^{\prime}=0.
        \end{cases}\nonumber
        \end{eqnarray}
  
  {\rm (10)} A Lagrangian $H$-umbilical submanifold of ratio $0$ and  non-constant mean 
 curvature 
 defined by $\pi\circ\psi$, where
   \bea
\begin{split}
\psi(s, u_2, \ldots, u_n)=
&\frac{\sqrt{\cosh{2s}}\exp[i{\tan}^{-1}(\tanh{s})]}{2}\biggl(1+\sum_{j=2}^nu_j^2+\sech 2s
-i\tanh 2s,\label{psi}\\
& i\sum_{j=2}^nu_j^2+i\sech 2s-i
+\tanh 2s, u_2, \ldots, u_n\biggr).\nonumber
\end{split}
\eea
  
  {\rm (11)} A Lagrangian $H$-umbilical submanifolds of ratio $(7-n)/3$ and  non-constant mean 
 curvature 
 defined by $\pi\circ\psi$, where $n\ne 7$ and
      \bea
\begin{split}
\psi(s, u_2, \ldots, u_n)=
&\cosh^{\frac{3}{n-1}}\biggl(\frac{n-1}{3}s\biggr)\exp\biggl[\frac{6i}{n-1}{\tan}^{-1}\biggl(\tanh\Bigl(\frac{n-1}{6}s\Bigr)\biggr)\biggr]\\
\times &\Biggl(\frac{1}{2}+\frac{1}{2}\sum_{j=1}^nu_j^2+\frac{1}{2}\sech^{\frac{6}{n-1}}
\biggl(\frac{n-1}{3}s\biggr)
-i\int_0^s\sech^{\frac{n+5}{n-1}}\biggl(\frac{n-1}{3}x\biggr)dx,\label{psi}\\
&\frac{i}{2}\sum_{j=2}^nu_j^2+
\frac{i}{2}\sech^{\frac{6}{n-1}}
\biggl(\frac{n-1}{3}s\biggr)-\frac{i}{2}
+\int_0^s\cosh^{\frac{n+5}{n-1}}\biggl(\frac{n-1}{3}x\biggr)dx, u_2, \ldots, u_n\Biggr).\nonumber
\end{split}
\eea

  {\rm (12)} Lagrangian $H$-umbilical submanifolds of non-constant mean 
 curvature 
 defined by $\pi\circ\psi$, where 
  \be
  \begin{split}
  \psi(s, u_2, \ldots, u_n)
  =&\exp\biggl[\int_0^s(k+i\mu)dt\biggr]\Biggl( 1+\dfrac{1}{2}\sum_{j=2}^nu_j^2-
  \int_0^s(k+i\mu){\rm exp}\biggl(-\int_0^x2kdt\biggr)dx,\nonumber\\
 & [-k(0)+i\mu(0)]\biggl[\dfrac{1}{2}\sum_{j=2}^nu_j^{\Add{2}}-\int_0^s(k+i\mu)\exp\biggl(-\int_0^x2kdt\biggr)dx\biggr], 
  u_2, \ldots, u_n\Biggr)
 \end{split} \nonumber
  \ee
 and $k(s)=\mu^{\prime}(s)/(\lambda(s)-2\mu(s))$,
  $\lambda\ne r\mu$ for any $r\in \mathbb{R}$, and moreover
    $\lambda=\lambda(s)$ and $\mu=\mu(s)$ are non-constant solutions of 
    \begin{eqnarray}
    \begin{cases}
    \mu^{\prime 2}=(\lambda-2\mu)^2(1-\mu^2), \\
     (3\lambda+(n-1)\mu)(\lambda-2\mu)\lambda^{\prime}+(n-1)\lambda
     (3\lambda+(n-5)\mu)\mu^{\prime}=0.
    \end{cases}\nonumber
    \end{eqnarray}
   \end{theorem}

 {\bf Proof of Theorems \ref{ct1}-\ref{ct3}}:
 Let $M^n$ be a Lagrangian $H$-umbilical submanifold, whose second fundamental form
 takes the form (\ref{lagh}), in a complex space form $\tilde M^2(4\epsilon)$
 with $\epsilon\in\{-1, 0, 1\}$. Assume that the mean curvature vector field
 is nowhere vanishing.
 
 We put $\omega_i^j(e_k)=\<\nabla_{e_k}e_i, e_j\>$. Then by using the equation
 (\ref{Codazzi}) of Codazzi, we have the following (cf. \cite{sa2}):
 \begin{align}
 & e_j\lambda=(2\mu-\lambda)\omega_j^1(e_1),\quad j>1, \label{c1}\\
 & e_1\mu=(\lambda-2\mu)\omega_1^2(e_2)=\cdots=(\lambda-2\mu)\omega_1^n(e_n),\label{c2}\\
 & e_j\mu=3\mu\omega_1^j(e_1),\quad j>1,\label{c3}\\
 & \mu\omega_1^j(e_1)=0, \quad j>1 \quad ({\rm for}\hskip5pt  n\geq 3),\label{c4}\\ 
 & \mu\omega_1^2(e_2)=\cdots=\mu\omega_1^n(e_n),\label{c5}\\
 & \mu\omega_1^i(e_j)=0, \quad i\ne j>1, \quad ({\rm for}\hskip5pt  n\geq 3), \label{c6}\\
 & (\lambda-2\mu)\omega_1^i(e_j)=0, \quad i\ne j>1, \quad ({\rm for}\hskip5pt  n\geq 3).\label{c7}
 \end{align}
 
It follows from  (\ref{tauH}) and (\ref{lagh}) that the tension field $\tau$ is given by
\be
\tau=(\lambda+(n-1)\mu)Je_1.\label{tau3}\ee
 For simplicity, we denote $A_{Je_j}$ by $A_j$.  
By using  (\ref{D}),  we have
 \be A_{D_{e_i}\tau} e_i=(e_i\lambda+(n-1)e_i\mu)A_1e_i+(\lambda+(n-1)\mu)
 \sum_{i=2}^n\omega_1^j(e_i)A_je_i. \label{AD}
 \ee
 The relation (\ref{relation}) implies that (\ref{lagh}) is equivalent to
  \be
   \begin{split}
   & A_1e_1=\lambda e_1,\\
   & A_1e_j=\mu e_j,\quad A_j e_1=\mu e_j, \quad A_je_j=\mu e_1, \quad j\geq 2,\label{A}\\
   & A_j e_k=0, \quad 2\leq j\ne k\leq n.
   \end{split}
   \ee 
  Combining (\ref{AD}) and (\ref{A}) we obtain 
 \be
 \begin{split}
 \sum_{i=1}^{n}A_{D_{e_i}\tau}e_i=&\sum_{i=1}^{n}(e_i\lambda+(n-1)e_i\mu)A_1e_i
 +(\lambda+(n-1)\mu)\sum_{i=1}^{n}\sum_{j=2}^n\omega_1^j(e_i)A_je_i \label{A_D}\\
 =&\sum_{i=1}^{n}(e_1\lambda+(n-1)e_1\mu)\lambda e_1
 +\sum_{i=2}^{n}(e_i\lambda+(n-1)e_i\mu)\mu e_i\\
 &+(\lambda+(n-1)\mu)\mu\biggl(\sum_{i=2}^n\omega_1^i(e_1)e_i+
  \sum_{i=2}^n\omega_1^i(e_i)e_1\biggr).
 \end{split}
 \ee
 We substitute  (\ref{tau3}) and (\ref{A_D}) into (\ref{TBE}). 
 Then,  by decomposing (\ref{TBE}) into the directions of $e_1$ and $e_j$ ($j>1$) we obtain 
 the following:
 \begin{align}
&(3\lambda+(n-1)\mu)(e_1\lambda+(n-1)e_1\mu)+2\mu
(\lambda+(n-1)\mu)\sum_{i=2}^n\omega_1^i(e_i)=0, \label{TBE1}\\
 &(\lambda+(n+1)\mu)(e_j\lambda+(n-1)e_j\mu)+2\mu
 (\lambda+(n-1)\mu)\omega_1^j(e_1)=0,
 \quad j>1.\label{TBE2}
 \end{align}

 {\bf Case (A):} $\mu=0$ and $\lambda\ne 0$. 
 
 Equations  (\ref{TBE1}) and (\ref{TBE2}) yield that $\lambda$ is constant. It follows from (\ref{c1}), (\ref{c2}) and (\ref{c7}) that
 $\omega_1^i(e_j)=0$ for all $i$, $j$. By the equation (\ref{Gauss}) of Gauss with 
 $X=W=e_i$ $(i\ne 1)$ and $Y=Z=e_1$, we have $\epsilon=0$. Hence we obtain case (1) of Theorem
 \ref{ct2}.

{\bf Case (B):} $\lambda=2\mu$ and $\mu\ne 0$.

By changing the sign of $e_1$ if necessary, we may assume that $\mu>0$.
By (\ref{c1}) and (\ref{c2})   we see that
 $\lambda$ is constant.  It follows from (\ref{c4})
  (\ref{c5}), (\ref{c6}),  (\ref{TBE1}) that $\omega_1^i(e_j)=0$ for all $i$, $j$. 
  In this case,  (\ref{TBE2}) is satisfied automatically.
 Similarly to the case of (A),   the equation (\ref{Gauss}) of Gauss implies that $\epsilon=-1$, $\mu=1$ and $\lambda=2$. Therefore, by the same argument as that in the proof of Theorem 6.1 in \cite{ch 3}
 we get case (1) of Theorem \ref{ct3}.

 {\bf Case (C):} $\mu(\lambda-2\mu)\ne 0$.
 
 Similary to  Case (B), we may assume that $\mu>0$.
For  $n=2$, substituting (\ref{c1}) and (\ref{c3}) into (\ref{TBE2}) leads to
 \be
 (\lambda+5\mu)\omega_1^2(e_1)=0.\label{lambda5}
 \ee
 Combining (\ref{c1}), (\ref{c3}) and  (\ref{lambda5}) yields
 $\omega_1^2(e_1)=e_2\lambda=e_2\mu=0$,
 i.e., the integral curves of $JH$ are geodesics in $M^2$. 
Thus we can apply Theorem \ref{t2}, \ref{t4} and \ref{t6}.  
 For $n\geq 3$,  by
 (\ref{c4}), (\ref{c1}) and (\ref{c3}) we have $\omega_1^i(e_1)=e_i\lambda=e_i\mu=0$ for all $i>1$.

Therefore, in both cases,
 (\ref{TBE2}) is satisfied automatically 
and hence the tangentially biharmonic condition
$(\ref{TBE})$
is equivalent to (\ref{TBE1}).
By using (\ref{c2}),  we see that (\ref{TBE1}) can be written as 
\be
(3\lambda+(n-1)\mu)(\lambda-2\mu)e_1\lambda+(n-1)\lambda
     (3\lambda+(n-5)\mu)e_1\mu=0.\label{TBE3}
\ee

We put $k=e_1\mu/(\lambda-2u)$. Then by combining (\ref{c2}) and the equation (\ref{Gauss}) of Gauss with 
 $X=W=e_i$ $(i\ne 1)$ and $Y=Z=e_1$, we have
\be
-e_1k-k^2=\epsilon+\lambda\mu-\mu^2.
\label{Gauss2}
\ee


{\bf Case (C.1):} $\mu$ is constant. 

In this case, $k=0$ and  hence (\ref{Gauss2}) leads to
 \be
 \epsilon+\lambda\mu-\mu^2=0, \label{gauss3}
 \ee
which implies that  $\lambda$ is also constant and 
$M^n$ has the ratio  $(1-\epsilon\mu^{-2})$. For $\epsilon=1$, applying Theorem \ref{t1} and \ref{t2}
gives case (1) of Theorem \ref{ct1}.
If  $\epsilon=-1$, then (\ref{gauss3}) and the condition $\lambda\ne 2\mu$ yield $\mu\ne 1$.
Therefore, by (1) and (2) of Theorem \ref{t5} and Theorem \ref{t6}, we get cases (2) and
(6) of Theorem \ref{ct3}.
For $\epsilon=0$, applying Theorem \ref{t3}, Theorem \ref{t4} and Remark \ref{remark2}, we obtain
case (2) of Theorem \ref{ct2}.

{\bf Case (C.2):} $\mu$ is non-constant and $\lambda=r\mu$ for some constant $r\ne  2$.

In this case, we have
  \be k=(\ln\hbox{\Erase{$|\mu|$}}\Add{\mu})^{\prime}/(r-2). \label{kk}\ee
 By (\ref{TBE3})  we get $r(r+n-1)(3r+n-7)=0$. Since the mean curvature vector field is nowhere
 vanishing, by (\ref{tau3}) we obtain 
 \be 
 r\in \{0, (7-n)/3\}.\label{ratio}
 \ee 


For $\epsilon=1$,  solving   (\ref{legen1}) for $\mu^{\prime}$ gives
\be
\mu^{\prime 2}=-(r-2)^2\mu^2(\mu^2+1)+c\mu^{\frac{2(r-3)}{r-2}}.\label{l1}
\ee 
It follows from (\ref{kk}) and (\ref{l1}) that
\be
(r-2)^2\mu^2(\mu^2+k^2+1)=c\mu^{\frac{2(r-3)}{r-2}},\nonumber 
\ee
which yields
$c>0$. Substituting (\ref{kk}) 
 into $(z_1, z_2)$ in Theorem \ref{t1}, using (\ref{ratio}) and (\ref{l1}),
 we get  cases (2) and (3) of Theorem \ref{ct1}.

 For  $\epsilon=-1$ and $\mu^2+k^2-1\ne 0$, solving  (\ref{legen2}) for $\mu^{\prime}$ gives
\be
\mu^{\prime 2}=(r-2)^2\mu^2(1-\mu^2)+c\mu^{\frac{2(r-3)}{r-2}}.\label{l2}
\ee 
This  implies that 
\be (r-2)^2\mu^2(\mu^2+k^2-1)=c\mu^{\frac{2(r-3)}{r-2}},\nonumber\ee and hence 
 the sign of $c$ coincides with the one of $\mu^2+k^2-1$.
 Substituting (\ref{kk}) 
 into $(z_1, z_2)$ in (1) and (2) of Theorem \ref{t5}, using (\ref{ratio}) and  (\ref{l2}), we obtain
cases (3) and (4)  of Theorem \ref{ct3} for $c>0$,   and cases (7) and (8)  of Theorem \ref{ct3} for $c<0$.

For $\epsilon=-1$ and $\mu^2+k^2-1=0$, by $(\ref{kk})$ we have
\be
\mu^{\prime 2}=(r-2)^2\mu^2(1-\mu^2).\label{muprime}
\ee
By solving (\ref{muprime})  and using (\ref{kk}), we get
 \be
 \mu={\rm sech}((r-2)s+c), \quad  k=-\tanh ((r-2)s+c)
 \ee
for some constant $c$.  
Substituting them  into $\psi$ of $(3)$ of Theorem \ref{t5},
after  suitable coordinate transformation, it becomes 
\bea
\begin{split}
\psi(s, u_2, \ldots, u_n)=
&\cosh^{-\frac{1}{r-2}}((r-2)s)\exp\biggl[\frac{2i}{r-2}{\tan}^{-1}\biggl(\tanh\Bigl(\frac{r}{2}-1\Bigr)s\biggr)\biggr]\\
\times &\biggl(\frac{1}{2}+\frac{1}{2}\sum_{j=1}^nu_j^2+\frac{1}{2}\cosh^{\frac{2}{r-2}}((2-r)s)
-i\int_0^s\cosh^{\frac{4-r}{r-2}}((r-2)x)dx,\label{psi}\\
&\frac{i}{2}\sum_{j=2}^nu_j^2+
\frac{i}{2}\cosh^{\frac{2}{r-2}}((2-r)s)-\frac{i}{2}
+\int_0^s\cosh^{\frac{4-r}{r-2}}((r-2)x)dx, u_2, \ldots, u_n\biggr).\nonumber
\end{split}
\eea
This expression and  (\ref{ratio})  give cases (10) and (11) of Theorem \ref{ct3}.

 For  $\epsilon=0$, by Theorem \ref{t3}, Theorem \ref{t4}, \Erase{Remark 4.2} and (\ref{ratio})
  we obtain cases (3) and (4) of Theorem \ref{ct2}.

{\bf Case (C.3):}   $\mu$ is non-constant and  $\lambda\ne r \mu$ for any $r\in\mathbb{R}$.

 For $\epsilon=1$, by 
 Theorem \ref{t1}, Theorem \ref{t2} and (\ref{TBE3}) we get case
 (4) of Theorem \ref{ct1}.
 If  $\epsilon=0$,  then  it follows from 
 $\lambda=\kappa$, $\mu=\theta^{\prime}$ and (\ref{theta2})
 that (\ref{TBE3}) is equivalent to (\ref{kappatheta}).
Hence, case (5) of Theorem \ref{ct2} is obtained.
For $\epsilon=-1$, by using Theorem \ref{t5}, Theorem \ref{t6} and (\ref{muprime}) we obtain
cases (5), (9) and (12) of Theorem \ref{ct3}.

The converse can be verified by a straightforward computation.
\qed\vspace{1.5ex}


A submanifold $M$ is said to be {\it isotropic} if at each point $p\in M$, $||h(v, v)||^2$ 
 is independent of the unit vector $v\in T_pM$ (see \cite{on}). 
Let $M^n$ be a non-minimal Lagrangian submanifold in a complex space form $\tilde M^n(4\epsilon)$
with $n\geq 3$. 
Then,  $M^n$ is isotropic if and only if $M^n$ is a Lagrangian $H$-umbilical submanifold of ratio 
$-1$ (see \cite{lw} and \cite{vr}).
 By applying Theorem \ref{ct1}-\ref{ct3}, we have
 \begin{corollary}
A complex space form $\tilde M^n(4\epsilon)$ with $n\geq 3$ 
admits isotropic tangentially biharmonic Lagrangian submanifolds with non-constant mean curvature
if and only if $n=10$.
 \end{corollary}


 \begin{remark}
{\rm As we have seen in Section $4$, non-minimal twistor holomorphic Lagrangian surfaces in $\mathbb{C}P^2$ and 
  the Whitney's sphere in $\mathbb{C}^n$ are Lagrangian $H$-umbilical submanifolds of ratio $3$. Thus Theorem \ref{ct1} and \ref{ct2} show that these submanifolds
  are not tangentially biharmonic.}
 \end{remark}

 \begin{remark}
 {\rm The submanifolds described in (1) of Theorem \ref{ct1} with 
 $$\mu=\sqrt{\frac{n+5\pm\sqrt{n^2+6n+25}}{2 n}}$$ are the only non-minimal
 biharmonic Lagrangian $H$-umbilical 
 submanifolds with constant mean curvature in complex space forms (cf. \cite[Theorem 7 and Remark 9]{sa2}).}
 \end{remark}

   \begin{remark}
  {\rm Each differential equation for $\lambda$ and $\mu$
  described in Theorem \ref{ct1} and \ref{ct3} can be transformed into an 
  autonomous system.
  Also, by (\ref{kappa}),  the differential equation described in \Erase{(6)} \Add{(5)} of Theorem \ref{ct2} can be rewritten as 
the third-order differential equation for $\alpha(s)$, and moreover, it  can be transformed  into an 
  autonomous system.
Therefore,  
by applying  Picard's existence  theorem for autonomous systems, we see that
there exist infinity many tangentially biharmonic Lagrangian $H$-umbilical submanifolds
of Theorem \ref{ct1}, Theorem \ref{ct3} and of type \Erase{(6)} \Add{(5)} of Theorem \ref{ct2}.}
   \end{remark}
 
 \begin{remark}
 {\rm By an argument  similar to that
  in \cite[Remark 1 \Erase{and 4}]{ch 4},
  a curve described in (4) \Erase{and (5)} of Theorem \ref{ct2} exists for any $n\geq 2$.}
 \end{remark}


\section{Corrections to this paper (added on June 18, 2023)}


\begin{itemize}
\setlength{\leftskip}{-0.5cm}
\item Theorems 4.1-4.6 should be replaced by respectively
\end{itemize}
 
 \begin{theorem}
 Let $f: M^n\rightarrow \mathbb{C}P^n(4)$ be a  Lagrangian $H$-umbilical immersion
 whose second fundamental form takes the form  $(\ref{lagh})$ with  $\mu(\lambda-2\mu)\ne 0$, where $n\geq 3$. 
 We put $k=e_1\mu/(\lambda-2\mu)$.
 Then, up
 to rigid motions of $\mathbb{C}P^n(4)$, $f$ is given by $\pi\circ \psi$, where
 $\pi: S^{2n+1}(1)\subset\mathbb{C}^{n+1}\rightarrow \mathbb{C}P^n(4): 
 w\mapsto w\cdot \mathbb{C}^{*}$ is the Hopf fibration and 
  $\psi$ is given by
 \be
 \psi(s, y_1, \ldots, y_n)=(z_1(s), z_2(s)y_1, \ldots, z_2(s)y_n), \quad y_1^2+\cdots+y_n^2=1,\nonumber
 \ee
 where $z=(z_1, z_2)$ is a unit speed  Legendre curve in $S^3(1)$
 given by
 \begin{align}
&z_1=\frac{i\mu(s)-k(s)}{\sqrt{\mu(s)^2+k(s)^2+1}}
\exp\biggl(\int_0^s(i\lambda(t)-i\mu(t))dt\biggr), \nonumber \\
&z_2= \frac{1}{\sqrt{\mu(s)^2+k(s)^2+1}}\exp\biggl(\int_0^si\mu(t)dt\biggr).
 \nonumber
 \end{align}
 Here $\mu=\mu(s)$ satisfies \be
  \mu^{\prime\prime}(\lambda-2\mu)-\mu^{\prime}(\lambda^{\prime}-3\mu^{\prime})
  +(\lambda-2\mu)^2(-\mu^2+\lambda\mu+1)=0.
  \ee
 \end{theorem}

 \begin{theorem}
 Let $f: M^2\rightarrow \mathbb{C}P^2(4)$ be a Lagrangian $H$-umbilical  immersion whose second fundamental form takes the form  $(\ref{lagh})$ with  $\mu(\lambda-2\mu)\ne 0$. 
 If the integral curves of $e_1$ are geodesics in $M^2$, then 
 up
  to rigid motions of $\mathbb{C}P^2(4)$, $f$ is given by the immersion described in Theorem $6.1$ with $n=2$.
 \end{theorem}

 \begin{theorem}
 Let $f: M^n\rightarrow \mathbb{C}^n$ be a Lagrangian $H$-umbilical  immersion whose second fundamental form takes the form  $(\ref{lagh})$ with $\mu(\lambda-2\mu)\ne 0$, where $n\geq 3$. 
  Then  up
  to rigid motions of $\mathbb{C}^n$, $f$ is   a complex extensor of the unit hypersphere in $\mathbb{E}^n$
  via a curve $F=F(s)$ in $\mathbb{C}^{*}$ whose curvature function $\kappa(s)$ and argument $\theta(s)$ satisfy
   $\kappa(s)=\lambda(s)$ and
     $\theta^{\prime}(s)=\mu(s)$. 
 \end{theorem}

 \begin{theorem}
  Let $f: M^2\rightarrow \mathbb{C}^2$ be a  Lagrangian $H$-umbilical immersion
  whose second fundamental form takes the form  $(\ref{lagh})$ with  $\mu(\lambda-2\mu)\ne 0$.
   If the  integral curves of $e_1$ are geodesics in $M^2$, then up
   to rigid motions of $\mathbb{C}^2$, $f$ is a complex extensor of the unit circle in $\mathbb
  {E}^2$
   via  a curve $F=F(s)$ in $\mathbb{C}^{*}$ whose
   curvature function $\kappa(s)$ and argument $\theta(s)$ satisfy
      $\kappa(s)=\lambda(s)$ and
        $\theta^{\prime}(s)=\mu(s)$. 
\end{theorem}

\begin{theorem}
 Let $f: M^n\rightarrow \mathbb{C}H^n(-4)$ be a  Lagrangian $H$-umbilical immersion
 whose second fundamental form takes the form  $(\ref{lagh})$ with $\mu(\lambda-2\mu)\ne 0$, where $n\geq 3$.
 We put $k=e_1\mu/(\lambda-2\mu)$.    
Then  up
 to rigid motions of $\mathbb{C}H^n(4)$, $f$ is  given by the composition $\pi\circ \psi$, where 
 $\pi: H_1^{2n+1}(-1)\subset \mathbb{C}^{n+1}_1\rightarrow \mathbb{C}H^n(-4):
 w\mapsto w\cdot \mathbb{C}^{*}$ is the Hopf fibration and 
  $\psi$ is one
 of the following immersion:
 
 {\rm (1)} $\mu^2+k^2-1>0$ and 
\be
\psi (s, y_1, \ldots, y_n)=(z_1(s), z_2(s)y_1, \ldots, z_2(s)y_n), \quad y_1^2+y_2^2+\cdots+y_n^2=1,\nonumber
 \ee
 where $z=(z_1, z_2)$ is a unit speed Legendre curve in $H_1^3(-1)$ given by
  \begin{align}
&z_1=\frac{i\mu(s)-k(s)}{\sqrt{\mu(s)^2+k(s)^2-1}}
\exp\biggl(\int_0^s(i\lambda(t)-i\mu(t))dt\biggr),\nonumber\\
&z_2= \frac{1}{\sqrt{\mu(s)^2+k(s)^2-1}}\exp\biggl(\int_0^si\mu(t)dt\biggr).
 \nonumber
 \end{align}
  
 {\rm (2)} $\mu^2+k^2-1<0$ and 
  \be
  \psi(s, y_1, \ldots, y_n)=(z_1(s)y_1, \ldots,  z_1(s)y_n, z_2(s)),  \quad y_1^2-y_2^2-\cdots-y_n^2=1,\nonumber
  \ee
  where $z=(z_1, z_2)$ is a unit speed Legendre curve in $H_1^3(-1)$ given by
   \begin{align}
&z_1= \frac{1}{\sqrt{1-\mu(s)^2-k(s)^2}}\exp\biggl(\int_0^si\mu(t)dt\biggr),\nonumber\\
&z_2=\frac{i\mu(s)-k(s)}{\sqrt{1-\mu(s)^2-k(s)^2}}
\exp\biggl(\int_0^s(i\lambda(t)-i\mu(t))dt\biggr).\end{align}

 {\rm(3)} $\mu^2+k^2-1=0$ and 
 \be
  \begin{split}
&\psi(s, u_2, \ldots, u_n)\\
=&\exp\biggl[\int_0^s(k(t)+i\mu(t))dt\biggr]\Biggl( 1+\dfrac{1}{2}\sum_{j=2}^nu_j^2-
  \int_0^s(k(x)+i\mu(x)){\rm exp}\biggl(-\int_0^x2k(t)dt\biggr)dx,\nonumber\\
 & [-k(0)+i\mu(0)]\biggl[\dfrac{1}{2}\sum_{j=2}^nu_j^2-\int_0^s(k(x)+i\mu(x))
 \exp\biggl(-\int_0^x2k(t)dt\biggr)dx\biggr], 
  u_2, \ldots, u_n\Biggr).
 \end{split} 
  \ee

  In the case of {\rm (1)} and {\rm (2)}, $\mu=\mu(s)$ satisfies
  \be
  \mu^{\prime\prime}(\lambda-2\mu)-\mu^{\prime}(\lambda^{\prime}-3\mu^{\prime})
  +(\lambda-2\mu)^2(-\mu^2+\lambda\mu-1)=0.
  \ee
 \end{theorem}

 \begin{theorem}
 Let $f: M^2\rightarrow \mathbb{C}H^2(-4)$ be a Lagrangian $H$-umbilical  immersion
 whose second fundamental form takes the form  $(\ref{lagh})$ with  $\mu(\lambda-2\mu)\ne 0$.
 If  the integral curves of $e_1$ are geodesics in $M^2$, then 
up
 to rigid motions of $\mathbb{C}H^2(-4)$,  $f$ is given by 
 one of the immersions described in Theorem $6.5$ with $n=2$.
 \end{theorem}

\begin{itemize}

\setlength{\leftskip}{-0.5cm}
\item In  Case (9) of Theorem 5.3, $z_1$ and $z_2$ should be replaced by respectively
 \begin{align}
 &z_1= \frac{1}{\sqrt{1-\mu^2-k^2}}\exp\biggl(\int_0^si\mu dt\biggr), \nonumber\\
&z_2=\frac{i\mu-k}{\sqrt{1-\mu^2-k^2}}
\exp\biggl(\int_0^s(i\lambda-i\mu)dt\biggr).
 \nonumber
 \end{align} 
\end{itemize}

 \end{document}